\documentclass[12pt]{article}
\usepackage{amsmath}
\usepackage{amssymb}
\renewcommand\smallskip{\vskip\smallskipamount}
\renewcommand\medskip{\vskip\medskipamount}
\renewcommand\bigskip{\vskip\bigskipamount}

\newcommand{\qed}{\hfill $\Box$ \medskip}

\begin{document}

\begin{center}
\begin{Large}
\textit{Counterexamples to the Local Solvability of
Monge-Amp\`{e}re Equations in the Plane}
\end{Large}
\bigskip

\begin{large}
\textsc{Marcus A. Khuri}
\end{large}

\end{center}

\medskip
\medskip

\indent\indent\indent\textsc{Abstract.}  In this paper, we present
$C^{\infty}$ examples of degenerate\newline\indent\indent\indent
hyperbolic and mixed type Monge-Amp\`{e}re equations in the plane,
\newline\indent\indent\indent which do not admit a local $C^{3}$
solution.

\medskip
\bigskip
\bigskip
\bigskip
\begin{center}
\textsc{0.  Introduction}
\end{center}

   Consider the class of two-dimensional Monge-Amp\`{e}re equations:
\begin{equation}
(u_{xx}+a(p,u,\nabla u))(u_{yy}+c(p,u,\nabla
u))-(u_{xy}+b(p,u,\nabla u))^{2}=f(p,u,\nabla u),
\end{equation}
where $p=(x,y)$.  The question of local solvability is to ask,
given smooth functions $a,b,c$, and $f$ defined in a neighborhood
of a point, say $(x,y)=0$, does there always exist a $C^{2}$
function $u(x,y)$, defined in a possibly smaller domain, which
satisfies (0.1)?  Note that we do not ask for $u(x,y)$ to satisfy
any boundary/initial conditions, have higher regularity, or to be
given in a predetermined domain. This is the most elementary
question that one can ask of a differential equation. Yet, it is
remarkable that the basic question of whether there exist any
examples of local nonsolvability, has remained open for this
well-studied class of equations.  The purpose of this paper is to
provide such examples.\par
   We first recall the known results.  Since (0.1) is elliptic if
$f>0$, hyperbolic if $f<0$, and of mixed type if $f$ changes sign,
the manner in which $f$ vanishes will play the primary role in the
hypotheses of any result.  The classical results state that a
solution always exists in the case that $f$ does not vanish at the
origin or is analytic ($a,b$, and $c$ are also required to be
analytic); these results follow easily from standard elliptic and
hyperbolic theory when $f$ does not vanish, and from the
Cauchy-Kovalevskaya theorem in the case that $a,b,c$ and $f$ are
analytic.  If $f(x,y,u,\nabla u)=f_{1}(x,y)f_{2}(x,y,u,\nabla u)$
with $f_{2}>0$, then C.-S. Lin provides an affirmative answer in
[7] and [8], when $f_{1}\geq 0$ or when $f_{1}(0)=0$ and $\nabla
f_{1}(0)\neq 0$.  When $f_{1}\leq 0$ and $\nabla f_{1}$ possesses
a certain nondegeneracy, Han, Hong, and Lin [3] show that a
solution always exists. Furthermore, in [4] and [5] the author
provides an affirmative answer in the case that $f_{1}$ has a
nondegenerate critical point at the origin, or degenerates to
arbitrary finite order along a single smooth curve (see also Han's
result [2]). Here, we shall prove\medskip

\textit{\textbf{Theorem.}  There exist sign changing and
nonpositive $f\in C^{\infty}(\mathbb{R}^{5})$, and $a,b,c\in
C^{\infty}(\mathbb{R}^{5})$, such that equation} (0.1)
\textit{possesses no $C^{3}$ solution in any neighborhood of the
origin.}\medskip

   The results mentioned above stem from work on a well-known
problem in geometry, namely the local isometric embedding problem
for two-dimensional Riemannian manifolds.  This problem is
equivalent to the local solvability of the following
Monge-Amp\`{e}re equation:
\begin{equation}
\det\nabla_{ij}u=K(\det g)(1-|\nabla_{g} u|^{2}),
\end{equation}
where $g$ is a given smooth Riemannian metric, K is its Gaussian
curvature, $\nabla_{ij}$ are second covariant derivatives, and
$\nabla_{g}$ is the gradient with respect to $g$.  Recently, N.
Nadirashvili and Y. Yuan have proposed counterexamples to the
isometric embedding problem in [9] and [10].  An immediate
consequence is the local nonsolvability of equation (0.2).
Although this observation was not mentioned by the authors, it is
quite significant since it represents the first nontrivial example
of a fully nonlinear equation exhibiting the property of local
nonsolvability.  (Of course in the setting of linear equations,
this phenomenon has received much attention through the work of
H\"{o}rmander, Nirenberg, Treves, and others since its original
discovery in 1957 by H. Lewy [6].)  The main distinction between
our theorem and the results of Nadirashvili and Yuan, besides the
difference in equations considered, is the fact that the proof
presented here is very elementary and does not rely on any
geometric significance that the equation may possess (as is
exhibited with equation (0.2)); as a result it is possible that
the methods presented here may be generalized to other
Monge-Amp\`{e}re equations.\par
   In the remainder of this section, we will partially construct
the functions $a,b,c$, and $f$ of the theorem, as well as reduce
the proof of this theorem to the problem of showing that certain
second derivatives of any solution of (0.1) must vanish along the
boundary of a sequence of squares.  Define sequences of disjoint
open squares $\{X^{n}\}_{n=1}^{\infty}$ and
$\{X^{n}_{1}\}^{\infty}_{n=1}$ whose sides are aligned with the
$x$ and $y$-axes, and such that $X^{n}$, $X^{n}_{1}$ are centered
at $q_{n}=(\frac{1}{n},0)$, $X^{n}\subset X^{n}_{1}$, and $X^{n}$,
$X^{n}_{1}$ have widths $\frac{1}{2n(n+1)}$, $\frac{1}{n(n+1)}$,
respectively. Set $a,b,c,f\equiv 0$ in $\mathbb{R}^{2}-\cup
X^{n}_{1}$. Define
\begin{equation*}
X=\{(x,y)\mid |x|<1,|y|<1\},
\end{equation*}
and let $\phi\in C^{\infty}(\overline{X})$ be such that $\phi$
vanishes to infinite order on $\partial X$, and either $\phi(p)>0$
or $\phi(p)<0$ for all $p\in X$; here $\overline{X}$ denotes the
closure of $X$.  We now define $a,b,c,f$ in $X^{n}$ by
$a,b,c\equiv 0$, and
\begin{equation}
f(p)=\gamma_{n}\phi(4n(n+1)(p-q_{n})),\text{ }\text{ }\text{
}\text{ }p\in \overline{X}^{\!\text{ }n},
\end{equation}
where $\{\gamma_{n}\}_{n=1}^{\infty}$ is a sequence of positive
numbers that will be chosen later, with the property that
$\lim_{n\rightarrow\infty}\gamma_{n}=0$.  In the next section, $f$
will be defined to be nonpositive in the remaining region
$\cup_{n=1}^{\infty}(X^{n}_{1}-X^{n})$.  Therefore, by choosing
$\phi$ to be positive or negative in each $X^{n}$, we obtain the
desired sign changing or nonpositive $f$ as mentioned in the
theorem.\par
   We now reduce the proof of the theorem as mentioned above.
Suppose that a local solution, $u\in C^{3}$, of (0.1) exists.  Let
$-\mathrm{v}_{n}$, $+\mathrm{v}_{n}$ represent the left and right
vertical portions of $\partial X^{n}$, respectively, and let
$+\mathrm{h}_{n}$, $-\mathrm{h}_{n}$ represent the top and bottom
horizontal portions of $\partial X^{n}$, respectively.  Now assume
that
\begin{equation}
u_{yy}|_{\pm\mathrm{v}_{n}}=0\text{ }\text{ }\text{ and }\text{
}\text{ } u_{xx}|_{\pm\mathrm{h}_{n}}=0, \text{ }\text{ }\text{
for all }\text{ }\text{ }n\geq N,
\end{equation}
where $N$ is the smallest integer such that $X^{N}$ is completely
contained within the domain of existence of $u$.  Let $n_{0}\geq
N$, and note that (0.3) and (0.4) imply that $u_{xy}|_{\pm
\mathrm{v}_{n_{0}}}=0$.  We may now integrate by parts to obtain a
contradiction,
\begin{equation*}
0\neq\int_{X^{n_{0}}}f=\int_{X^{n_{0}}}u_{xx}u_{yy}-u_{xy}^{2}
=\int_{\partial X^{n_{0}}}u_{xx}u_{y}\nu_{2}-u_{xy}u_{y}\nu_{1}=0,
\end{equation*}
where $(\nu_{1},\nu_{2})$ are the components of the unit outward
normal to $\partial X^{n_{0}}$.  Thus, our theorem is reduced to
the proof of (0.4).\par
   The outline of the paper is as follows.  In section $\S 1$ we
complete the construction of $a,b,c$, and $f$.  Furthermore,
assuming that (0.4) does not hold, we find a certain integral
equality that $u$ must satisfy.  In order to violate this integral
equality, we construct approximate solutions to a homogeneous
degenerate hyperbolic equation in section $\S 2$.
\begin{center}
\textsc{1.  The Integral Equality}
\end{center} \setcounter{equation}{0}
\setcounter{section}{1}

   The purpose of this section is to construct a sequence of
integral equalities, valid for $C^{3}$ solutions of (0.1) in
subdomains of $X_{1}^{n}-X^{n}$ if (0.4) is violated.  However,
before obtaining the integral equalities we will first complete
the construction of $a,b,c$, and $f$ in the regions
$X_{1}^{n}-X^{n}$. Extend the line segments $\pm\mathrm{v}_{n}$,
$\pm\mathrm{h}_{n}$ until they reach $\partial X_{1}^{n}$, so that
we obtain four rectangles each bounded by $\partial X^{n}$,
$\partial X_{1}^{n}$, and the extended segments
$\pm\mathrm{v}_{n}$, $\pm\mathrm{h}_{n}$. Denote the rectangles to
the left and right of $\partial X^{n}$ by $-V_{n}$, $+V_{n}$
respectively, and denote the rectangles to the top and bottom of
$\partial X^{n}$ by $+H_{n}$, $-H_{n}$ respectively.  We then set
$a,b,c,f\equiv 0$ in $(X_{1}^{n}-X^{n})-(\pm V_{n}\cup\pm H_{n})$,
define $a=a_{n}(x,y)uu_{y}$, $b=c\equiv 0$ in $\pm V_{n}$, and
$c=c_{n}(x,y)uu_{x}$, $a=b\equiv 0$ in $\pm H_{n}$, for some
$a_{n}\in C^{\infty}(\pm\overline{V}_{n})$ and $c_{n}\in
C^{\infty}(\pm\overline{H}_{n})$ to be given below.  Lastly, in
$\pm V_{n}\cup\pm H_{n}$ we will write
$f=K_{n}(x,y)+g_{n}(x,y,\nabla u)$ for nonpositive functions
$K_{n}\in C^{\infty}(\overline{\pm V_{n}\cup\pm H_{n}})$,
$g_{n}\in C^{\infty}(\overline{\pm V_{n}\cup\pm
H_{n}}\times\mathbb{R}^{2})$ also to be given below. In order to
motivate the construction of $f$ in the regions $\pm V_{n}\cup\pm
H_{n}$, we will now convert (0.1) into a quasilinear equation by
applying an appropriate Legendre transformation.\par
   Let $u\in C^{3}$ be a local solution of (0.1), and as above let
$n_{0}$ be such that $X^{n_{0}}_{1}$ is contained within the
domain of existence of $u$.  Let $p=(p_{1},p_{2})\in
+\mathrm{v}_{n_{0}}$, and assume that (0.4) is violated, so that
$u_{yy}(p)\neq 0$.  Then we have a well-defined $C^{2}$ Legendre
transformation $T:(x,y)\mapsto(\alpha,\beta)$ defined in a
sufficiently small neighborhood, $B_{p}$ of $p$, and given by
\begin{equation*}
\alpha=x-p_{1},\text{ }\text{ }\text{ }\text{ }\text{
}\beta=u_{y}(x,y).
\end{equation*}
It follows that $u$ must satisfy a quasilinear equation in the new
variables.\medskip

\textit{\textbf{Lemma 1.}  There exist constants $\alpha_{0}>0$,
$\beta_{1}>\beta_{*}>\beta_{0}$, and a rectangle
$D=(0,\alpha_{0})\times(\beta_{0},\beta_{1})\subset T(B_{p})$
where $T(p)=(0,\beta_{*})$, such that in $D$ we have
$|u_{y}u_{yy}|>0$ and}
\begin{equation}
Lu:=u_{\alpha\alpha}+(K_{n_{0}}u_{\beta})_{\beta}-(2\beta^{-1}K_{n_{0}}
+\beta^{2}a_{n_{0}}+\beta
uu_{yy}a_{n_{0}\beta})u_{\beta}=G_{n_{0}},
\end{equation}
\textit{where}
\begin{equation*}
G_{n_{0}}=-(g_{n_{0}}u_{\beta})_{\beta}
+2\beta^{-1}g_{n_{0}}u_{\beta}.
\end{equation*}

\textit{Proof.}  If $u_{y}(p)=0$, then we could instead take any
point $\overline{p}\in +\mathrm{v}_{n_{0}}$ near $p$ with
$u_{y}(\overline{p})\neq 0$ and $u_{yy}(\overline{p})\neq 0$.  The
existence of a rectangle $D$ in which $|u_{y}u_{yy}|>0$ now
follows.  Moreover, if we set
\begin{equation*}
F=f-au_{yy}+2bu_{xy}-cu_{xx}+b^{2}-ac,
\end{equation*}
then for any $\psi\in C_{c}^{\infty}(D)$ we have
\begin{eqnarray*}
\int_{D}\!(u_{\alpha}\psi_{\alpha}+\!Fu_{\beta}\psi_{\beta})d\alpha
d\beta\!\!\!\!&=&\!\!\!\!\!\int_{T^{-1}(D)}[(u_{x}\!-\!\frac{u_{xy}}{u_{yy}}u_{y})(\psi_{x}\!
-\!\frac{u_{xy}}{u_{yy}}\psi_{y})+\!F(\frac{1}{u_{yy}})^{2}u_{y}\psi_{y}]u_{yy}dxdy\\
&=&\!\!\!\!\int_{T^{-1}(D)}[(u_{x}u_{yy}-u_{xy}u_{y})\psi_{x}+(u_{y}u_{xx}
-u_{x}u_{xy})\psi_{y}]dxdy\\
&=&\!\!\!\!-2\int_{T^{-1}(D)}F\psi dxdy\\
&=&\!\!\!\!-2\int_{D}\beta^{-1}Fu_{\beta}\psi d\alpha d\beta.
\end{eqnarray*}
Recalling that $f=K_{n_{0}}+g_{n_{0}}$ and $b=c\equiv 0$ in
$+V_{n_{0}}$ we obtain
\begin{equation*}
u_{\alpha\alpha}+(K_{n_{0}}u_{\beta})_{\beta}-(2\beta^{-1}K_{n_{0}}+
(au_{yy}u_{\beta})_{\beta}u_{\beta}^{-1}-2\beta^{-1}au_{yy})u_{\beta}=G_{n_{0}},
\end{equation*}
from which (1.1) follows with $a=a_{n_{0}}uu_{y}$.\qed
   We may view (1.1) as a linear equation which possesses a
solution $u\in C^{2}(D)$.  Since $K_{n_{0}}$ is nonpositive in
$+V_{n_{0}}$, equation (1.1) is degenerate hyperbolic in D.  Let
$\overline{a},\overline{b},\overline{c},\overline{d},\overline{g}\in
C^{\infty}(D)$ with $\overline{a}\geq 0$, and consider the linear
degenerate hyperbolic equation:
\begin{equation}
z_{\alpha\alpha}-(\overline{a}z_{\beta})_{\beta}+\overline{b}z_{\beta}+
\overline{c}z_{\alpha}+\overline{d}z=\overline{g}.
\end{equation}
The local solvability of (1.2) is highly dependent upon certain
relationships between the coefficients $\overline{a}$ and
$\overline{b}$, the so called Levi conditions.  One of the most
powerful Levi conditions was given by Oleinik [11], who proved
that the Cauchy problem for (1.2), with data prescribed on the
line $\alpha=0$, is well-posed if there exists an integer $J>0$
and constants $A,B>0$,
$\alpha_{0}=0<\alpha_{1}<\cdots<\alpha_{J}$, such that either
\begin{equation}
B(\alpha-\alpha_{j-1})\overline{b}^{2}\leq
A\overline{a}+\overline{a}_{\alpha}\text{ }\text{ }\text{ }\text{
or }\text{ }\text{ }\text{ }
B(\alpha_{j}-\alpha)\overline{b}^{2}\leq[A+\frac{1}{B(\alpha_{j}-\alpha)}]
\overline{a}-\overline{a}_{\alpha}
\end{equation}
holds for $\alpha_{j-1}\leq\alpha\leq\alpha_{j}$, $j=1,\ldots,J$.
Note that (1.3) implies that either
\begin{equation*}
0\leq A\overline{a}+\overline{a}_{\alpha}\text{ }\text{ }\text{
}\text{ or }\text{ }\text{ }\text{
}0\leq[A+\frac{1}{B(\alpha_{j}-\alpha)}]\overline{a}-\overline{a}_{\alpha},
\end{equation*}
so that we must have either
$\overline{a}(\alpha,\overline{\beta})=0$ for all
$\alpha_{j-1}\leq\alpha\leq\overline{\alpha}$, or all
$\overline{\alpha}\leq\alpha\leq\alpha_{j}$, respectively, if
$\overline{a}(\overline{\alpha},\overline{\beta})=0$.  Therefore,
Oleinik's result does not allow the coefficient $\overline{a}$ to
have an infinite sequence of isolated zeros as $\alpha\rightarrow
0$, as would be the case if $\overline{a}$ exhibited fast
oscillations near $\alpha=0$.  In fact, counterexamples to the
local solvability of (1.2) have been found [1] in the case that
$\overline{a}$ possesses this type of behavior.\par
   With this intuition, we will define $K_{n}$ in $+V_{n}$ to have
special fast oscillations as $x\rightarrow\partial X^{n}$.  Let
$x_{n}=x-(\frac{1}{n}+\frac{1}{4n(n+1)})$, so that in the new
coordinates $+V_{n}$ is given by
\begin{equation*}
+V_{n}=\{(x_{n},y)\mid 0<x_{n}<\frac{1}{4n(n+1)},\text{
}\frac{-1}{4n(n+1)}<y<\frac{1}{4n(n+1)}\}.
\end{equation*}
Let $k_{n}$ be the smallest integer such that
$k_{n}>\frac{4n(n+1)}{\pi}$, and set
$I_{k}=(\frac{1}{\pi(k+1)},\frac{1}{\pi k})$,
$k\in\mathbb{Z}_{>0}$.  Then define a smooth function $K$ for
$x_{n}>0$ by
\begin{equation*}
K(x_{n})=\begin{cases} e^{-m_{k}^{-2}-\sin^{-2}(\frac{1}{x_{n}})}
& \text{if $x_{n}\in I_{k}$, $k\in
\mathbb{Z}_{k\geq k_{n}}$},\\
0 & \text{if $x_{n}\in\partial I_{k}$},\\
0 & \text{if $x_{n}\geq\frac{1}{\pi k_{n}}$},
\end{cases}
\end{equation*}
where $m_{k}$ denotes the unique zero of $\cos(\frac{1}{x_{n}})$
in the interval $I_{k}$.  The function $K$ will be used to give
$K_{n}$ fast oscillations, however $K_{n}$ will also be required
to have specific behavior in the $y$-direction as well.  In order
to accomplish this we partition $+V_{n}$ into small rectangles
\begin{equation*}
\mathcal{R}_{k,i,n}=\{(x_{n},y_{n})\mid x_{n}\in I_{k},\text{ }
ik^{-1/2}<y_{n}<(i+1)k^{-1/2}\},\text{ }\text{ }\text{ }\text{
}i=0,1,\ldots,\overline{i}(k,n),
\end{equation*}
where $y_{n}=y+\frac{1}{4n(n+1)}$ and $\overline{i}(k,n)$ denotes
the largest integer such that
$\mathcal{R}_{k,\overline{i}(k,n),n}\cap+V_{n}\neq\emptyset$. Next
let $\psi\in C^{\infty}(-\infty,\infty)$ be a nonnegative
$1$-periodic function such that
\begin{equation*}
\psi(\overline{y})=\begin{cases} 1 & \text{if
$\frac{1}{4}\leq\overline{y}\leq\frac{3}{4}$},\\
0 & \text{if $0\leq\overline{y}\leq\frac{1}{8}$ or
$\frac{7}{8}\leq\overline{y}\leq 1$},
\end{cases}
\end{equation*}
and set
\begin{equation*}
\Psi_{n}(x_{n},y_{n})=\begin{cases} \psi(k^{1/2}y_{n}) & \text{if
$(x_{n},y_{n})\in\mathcal{R}_{k,i,n}$, $0\leq
i<\overline{i}(k,n)$},\\
0 & \text{if
$(x_{n},y_{n})\in\mathcal{R}_{k,\overline{i}(k,n),n}$}.
\end{cases}
\end{equation*}
We then define
\begin{equation*}
K_{n}(x,y)=-\gamma_{n}K(x_{n})\Psi_{n}(x_{n},y_{n}),\text{ }\text{
}\text{ }\text{ }\text{ }(x,y)\in+V_{n},
\end{equation*}
where $\gamma_{n}$ was given in (0.3).  Note that $K_{n}\leq 0$ in
$+V_{n}$ and $K_{n}\in C^{\infty}(+\overline{V}_{n})$.
Furthermore, we define $K_{n}$ analogously in the rectangles
$-V_{n}$, $\pm H_{n}$, where in $\pm H_{n}$ the roles of $x$ and
$y$ are reversed.\par
   We now choose $a_{n}$, $c_{n}$, and $g_{n}$.  Set
$\xi_{1}(x_{n})=\sin^{-4}(\frac{1}{x_{n}})$ and define
\begin{equation*}
a_{n}(x,y)=2\xi_{1}^{'}(x_{n})\sqrt{\gamma_{n}K(x_{n})}\Psi_{n}(x_{n},y_{n}),
\text{ }\text{ }\text{ }\text{ }\text{ }(x,y)\in +V_{n}.
\end{equation*}
Note that $a_{n}\in C^{\infty}(+\overline{V}_{n})$.  We also
define $a_{n}$ in $-V_{n}$ and $c_{n}$ in $\pm H_{n}$ similarly,
where in $\pm H_{n}$ the roles of $x$ and $y$ are reversed. Next
let $\{N_{k}\}_{k=1}^{\infty}$ be a sequence of positive integers,
and let $\psi_{N_{k}}\in C^{\infty}(-\infty,\infty)$ be a
nonnegative $1$-periodic function such that
\begin{equation}
\int_{0}^{1}\psi_{N_{k}}^{'}(\overline{y})\psi^{(N_{k})}(\overline{y})
d\overline{y}\geq 1,\text{ }\text{ }\text{ }\text{ }\text{
}k\in\mathbb{Z}_{>0};
\end{equation}
here and below we use the notation
$\psi^{(j)}(\overline{y})=\frac{d^{j}}{d\overline{y}^{j}}\psi(\overline{y})$.
The integers $N_{k}$ are to be chosen converging to infinity
sufficiently slow so that
\begin{equation}
\lim_{k\rightarrow\infty}e^{-\frac{1}{2}m_{k}^{-2}}\sup_{\overline{y}\in
(0,1)}|\psi_{N_{k}}^{(j)}(\overline{y})|= 0,\text{ }\text{ }\text{
}\text{ for each }\text{ }\text{ }j\in\mathbb{Z}_{\geq 0}.
\end{equation}
If we did not require $\psi_{N_{k}}$ to be nonnegative, that is,
if in the theorem we only wished to construct examples of mixed
type Monge-Amp\`{e}re equations (see the definition of $g_{n}$
below), then we could simply take $\psi_{N_{k}}=\psi^{(N_{k}-1)}$,
since the Poincar\'{e} inequality implies that
\begin{equation*}
\int_{0}^{1}\psi^{2}(\overline{s})d\overline{s}\leq\pi^{-2N_{k}}
\int_{0}^{1}[\psi^{(N_{k})}]^{2}(\overline{s})d\overline{s},
\end{equation*}
where we have used the fact that $\pi^{2}$ is the principal
eigenvalue of $-\frac{d^{2}}{d\overline{s}^{2}}$ on the interval
$(0,1)$.  Let $\{\tau_{k}\}_{k=1}^{\infty}$ be a sequence of
positive numbers to be chosen later with
$\lim_{k\rightarrow\infty} \tau_{k}=\infty$, then for $(x,y)\in
+V_{n}$ we define
\begin{equation*}
g_{n}(x,y,\nabla u)=\begin{cases}
\tau_{k}^{-N_{k}}K_{n}(x,y)\psi_{N_{k}}(\tau_{k}s_{n}(x,y,u_{y}))
& \text{if
$(x_{n},y_{n})\in\mathcal{R}_{k,i,n}$, $0\leq i<\overline{i}(k,n)$},\\
0 & \text{if
$(x_{n},y_{n})\in\mathcal{R}_{k,\overline{i}(k,n),n}$},
\end{cases}
\end{equation*}
where $s_{n}=u_{y}(x,y)+\int_{0}^{x_{n}}
\sqrt{\gamma_{n}K(\overline{x})}d\overline{x}$.  Note that
$g_{n}\leq 0$ and according to (1.5), $g_{n}\in
C^{\infty}(+\overline{V}_{n}\times\mathbb{R}^{2})$.  In the
regions $-V_{n}$, $\pm H_{n}$, $g_{n}$ is defined similarly, where
again the roles of $x$ and $y$ are reversed in $\pm H_{n}$.
Lastly, by choosing a sequence $\{\gamma_{n}\}_{n=1}^{\infty}$
which converges to zero sufficiently fast we have $a,b,c,f\in
C^{\infty}(\mathbb{R}^{5})$.\par
   Now that $a,b,c$, and $f$ are fully defined, we shall obtain a
sequence of integral equalities for $u$ inside $+V_{n_{0}}$.  In
order to construct the domains in which these integral equalities
will be valid, we need the following change of coordinates
$T_{1}:(\alpha,\beta)\mapsto(t,s)$ given by
\begin{equation*}
t=\alpha,\text{ }\text{ }\text{ }\text{ }\text{ }s(\alpha,\beta)=
\beta+\int_{0}^{\alpha}\sqrt{\gamma_{n_{0}}K(\overline{\alpha})}d\overline{\alpha},
\end{equation*}
which is well-defined near $T(p)$.  Let
$\{\mathcal{R}_{k,i_{k},n_{0}}\}_{k=1}^{\infty}$ be a sequence of
rectangles such that
$\mathcal{R}_{k,i_{k},n_{0}}\cap\{(x,y)\in+V_{n_{0}}\mid
y=p_{2}\}\neq\emptyset$, so that this sequence converges to $p$.
Then for $(x,y)\in\mathcal{R}_{k,i_{k},n_{0}}$ we have
$|x-p_{1}|+|y-p_{2}|^{2}=O(k^{-1})$, which implies that
\begin{equation*}
s=u_{yy}(p)y+(u_{y}(p)-u_{yy}(p)p_{2})+O(k^{-1}).
\end{equation*}
Therefore if $\tau_{k}\geq k$ and $k$ is large enough, there exist
integers $J_{k}$ such that the rectangles in the $ts$-plane,
\begin{equation*}
R_{k,\tau_{k}}=\{(t,s)\mid t\in I_{k},\text{
}J_{k}<\tau_{k}s<J_{k}+1\},
\end{equation*}
satisfy
\begin{equation*}
T^{-1}(T_{1}^{-1}(R_{k,\tau_{k}}))\subset\{(x,y)\in+V_{n_{0}}\mid
\Psi_{n_{0}}(x_{n_{0}},y_{n_{0}})\equiv
1\}\cap\mathcal{R}_{k,i_{k},n_{0}}.
\end{equation*}
It follows that
$a_{n_{0}\beta}(x,y)=(a_{n_{0}y}u_{yy}^{-1})(x,y)=0$ for $(x,y)\in
T^{-1}(T_{1}^{-1}(R_{k,\tau_{k}}))$.  Then in $R_{k,\tau_{k}}$ the
expression for the operator $L$ of (1.1) becomes
\begin{equation*}
Lu=u_{tt}+2Au_{ts}+Bu_{s},
\end{equation*}
where the coefficients $A$ and $B$ are smooth in $ts$-coordinates
and are given by
\begin{equation*}
A=\sqrt{\gamma_{n_{0}}K},\text{ }\text{ }\text{ }\text{ }\text{ }
B=2\beta^{-1}\gamma_{n_{0}}K-2\beta^{2}\xi_{1}^{'}\sqrt{\gamma_{n_{0}}K}
+(\sqrt{\gamma_{n_{0}}K})_{t}.
\end{equation*}
Let $z_{k}\in C^{\infty}(\overline{R}_{k,\tau_{k}})$ vanish
identically on the boundary, then multiplying equation (1.1)
through by $z_{k}$ and integrating by parts gives the desired
integral equality
\begin{equation}
\int_{R_{k,\tau_{k}}}uL^{*}z_{k}dtds=\int_{R_{k,\tau_{k}}}G_{n_{0}}z_{k}dtds,
\end{equation}
where $L^{*}$ denotes the formal adjoint of $L$.
\begin{center}
\textsc{2.  Violation of the Integral Equality}
\end{center} \setcounter{equation}{0}
\setcounter{section}{2}

   In this section we will violate the integral equality (1.6)
for large $k$, by choosing the sequence
$\{\tau_{k}\}_{k=1}^{\infty}$ to grow sufficiently fast, and by
constructing approximate solutions of the homogeneous adjoint
equation $L^{*}z=0$ vanishing identically on $\partial
R_{k,\tau_{k}}$, so that the left-hand side of (1.6) tends to zero
much faster than the right-hand side.\par
   The approximate solutions will be of the form
\begin{equation*}
z_{k}=e^{\eta(t,s)}\sum_{i=0}^{N_{k}}\tau_{k}^{-i}a_{i}(t,s)b_{i}(\tau_{k}s):=e^{\eta}
\overline{z},
\end{equation*}
where $a_{i},\eta\in C^{\infty}(R_{k,\tau_{k}})$ and
$b_{i}=\psi^{(N_{k}-i)}$.  A calculation shows that
\begin{equation}
e^{-\eta}L^{*}z_{k}=\overline{z}_{tt}+2A\overline{z}_{ts}
+\overline{B}\overline{z}_{s}+\overline{C}\overline{z}_{t}+
\overline{D}\overline{z},
\end{equation}
where
\begin{equation*}
\overline{B}=2A_{t}-B+2A\eta_{t},\text{ }\text{ }\text{ }
\overline{C}=2(A_{s}+A\eta_{s}+\eta_{t}),\text{ }\text{ }\text{ }
\overline{D}=\eta_{t}^{2}+\eta_{tt}+\overline{B}\eta_{s}+\overline{B}_{s}.
\end{equation*}
In order to simplify (2.1) we choose $\eta$ so that
$\overline{B}\equiv 0$, that is we set
\begin{equation*}
\eta(t,s)=\int_{m_{k}}^{t}(-\beta^{2}\xi_{1}^{'}(\overline{t})
+\beta^{-1}\sqrt{\gamma_{n_{0}}K(\overline{t})})d\overline{t}+\frac{1}{4}\xi_{2}(t),
\end{equation*}
where $\xi_{1}(t)=\sin^{-4}(\frac{1}{t})$ and
$\xi_{2}(t)=\sin^{-2}(\frac{1}{t})$.  It follows that
\begin{equation}
C^{-1}e^{-\sigma_{1}\sin^{-4}(\frac{1}{t})}\leq e^{\eta} \leq
Ce^{-\sigma_{2}\sin^{-4}(\frac{1}{t})}\text{ }\text{ }\text{ in
}\text{ }\text{ }R_{k,\tau_{k}},
\end{equation}
for some constants $C,\sigma_{1},\sigma_{2}>0$ independent of $k$.
Furthermore, (2.1) becomes
\begin{eqnarray*}
e^{-\eta}L^{*}z_{k}&=&2\tau_{k}Aa_{0t}b_{0}^{'}\\
& &+\sum_{i=0}^{N_{k}-1}\tau_{k}^{-i}(a_{itt}+2Aa_{its}
+\overline{C}a_{it}+\overline{D}a_{i}+2Aa_{(i+1)t})b_{i}\\
&
&+\tau_{k}^{-N_{k}}(a_{N_{k}tt}+2Aa_{N_{k}ts}+\overline{C}a_{N_{k}t}
+\overline{D}a_{N_{k}})b_{N_{k}}.
\end{eqnarray*}
This suggests that we inductively choose
$a_{0}=\mathrm{sgn}(u_{y}u_{yy}^{-1}(p)):=\pm 1$,
\begin{equation*}
2Aa_{(i+1)t}=-a_{itt}
-2Aa_{its}-\overline{C}a_{it}-\overline{D}a_{i},\text{ }\text{
}\text{ }a_{i+1}(m_{k},s)=0,\text{ }\text{ }\text{
}i=0,\ldots,N_{k}-1.
\end{equation*}
Therefore (2.2) shows that $z_{k}\in
C^{\infty}(\overline{R}_{k,\tau_{k}})$ vanishes identically on
$\partial R_{k,\tau_{k}}$, and
\begin{equation}
\max_{j=0,1,2}\sup_{(t,s)\in
R_{k,\tau_{k}}}e^{\eta(t,s)}|\nabla^{j}a_{i}(t,s)|
(1+|\overline{C}(t,s)|+|\overline{D}(t,s)|)\leq M(i,k,u)
\end{equation}
for some constants $M(i,k,u)$ which only depend on $\tau_{k}$
through their dependence on $u$.  As a result, it is clear that
$\tau_{k}$ can be chosen so that for any $w\in
C^{3}(+\overline{V}_{n_{0}})$ with $w_{y}(p)\neq 0$ and
$w_{yy}(p)\neq 0$,
\begin{equation}
C_{1}(\psi_{N_{k}})\sum_{i=1}^{N_{k}}\tau_{k}^{-i}C_{N_{k}-i}(\psi)M(i,k,w)\leq
\overline{M}_{0}(w)e^{-k^{4}},
\end{equation}
\begin{equation}
\tau_{k}^{-1}M(N_{k},k,w)\leq\overline{M}_{1}(w)e^{-k^{4}},\text{
}\text{ }\text{ }\text{ }\text{
}\tau_{k}^{-1}C_{N_{k}}(\psi)C_{0}(\psi_{N_{k}})\leq k^{-4},
\end{equation}
where
\begin{equation*}
C_{j}(\psi):=\sup_{\overline{s}\in(0,1)}|\psi^{(j)}(\overline{s})|,\text{
}\text{ }\text{ }\text{ }\text{
}C_{j}(\psi_{N_{k}}):=\sup_{\overline{s}\in(0,1)}|\psi^{(j)}_{N_{k}}(\overline{s})|,
\end{equation*}
and the constants $\overline{M}_{0}(w)$ and $\overline{M}_{1}(w)$
depend only on $\beta_{0}(w)$ and $\beta_{1}(w)$ (see lemma 1).
Then combining (2.3) and (2.5) we have
\begin{equation}
|\int_{R_{k,\tau_{k}}}uL^{*}z_{k}dtds|\leq C\tau_{k}^{-N_{k}-1}
M(N_{k},k,u)\leq C\overline{M}_{1}(u)e^{-k^{4}}\tau_{k}^{-N_{k}}.
\end{equation}\par
   We now estimate the right-hand side of the integral equality
(1.6).  First observe that lemma 1 together with the definition of
$g_{n_{0}}$ yields
\begin{equation*}
G_{n_{0}}(t,s)=\gamma_{n_{0}}\tau_{k}^{-N_{k}+1}\psi_{N_{k}}^{'}(\tau_{k}s)
K(t)u_{y}u_{yy}^{-1}
+O(\tau_{k}^{-N_{k}}C_{0}(\psi_{N_{k}})e^{-m_{k}^{-2}}),\text{
}\text{ }\text{ }\text{ }(t,s)\in R_{k,\tau_{k}}.
\end{equation*}
Moreover, using (2.3) we have
\begin{equation*}
z_{k}(t,s)=\mathrm{sgn}(u_{y}u_{yy}^{-1}(p))e^{\eta}\psi^{(N_{k})}(\tau_{k}s)+
O(\sum_{i=1}^{N_{k}}\tau_{k}^{-i}C_{N_{k}-i}(\psi)M(i,k,u)).
\end{equation*}
We now apply (1.4), (2.2), (2.4), and (2.5) to obtain
\begin{eqnarray}
\int_{R_{k,\tau}}G_{n_{0}}z_{k}dtds\!&\geq&\!Ce^{-m_{k}^{-2}}
\tau_{k}^{-N_{k}}(\int_{\frac{1}{\pi(k+1)}}^{\frac{1}{\pi
k}}e^{-\sigma_{1}\sin^{-4}(\frac{1}{t})-\sin^{-2}(\frac{1}{t})}dt)
(\int_{0}^{1}\psi_{N_{k}}^{'}(\overline{s})
\psi^{(N_{k})}(\overline{s})d\overline{s})\nonumber\\
&
&\!+O(\tau_{k}^{-N_{k}-1}e^{-m_{k}^{-2}}C_{N_{k}}(\psi)C_{0}(\psi_{N_{k}})\\
& &\!+
\tau_{k}^{-N_{k}}e^{-m_{k}^{-2}}C_{1}(\psi_{N_{k}})\sum_{i=1}^{N_{k}}
\tau_{k}^{-i}C_{N_{k}-i}(\psi)M(i,k,u))\nonumber\\
&\geq&\! C_{1}k^{-2}e^{-m_{k}^{-2}}\tau_{k}^{-N_{k}}
+O(\tau_{k}^{-N_{k}}e^{-m_{k}^{-2}}\overline{M}_{0}(u)e^{-k^{4}})\nonumber,
\end{eqnarray}
for some constants $C,C_{1}>0$ independent of $k$.\par
   We may now complete the proof of the theorem.  Combining (1.6),
(2.6), and (2.7) produces
\begin{equation*}
C^{-1}k^{-2}e^{-m_{k}^{-2}}\tau_{k}^{-N_{k}}\leq
\int_{R_{k,\tau_{k}}}G_{n_{0}}z_{k}dtds=\int_{R_{k,\tau_{k}}}uL^{*}z_{k}dtds\leq
C\overline{M}_{1}(u)e^{-k^{4}}\tau_{k}^{-N_{k}},
\end{equation*}
which leads to a contradiction for large $k$.  It follows that
$u_{yy}(p)=0$.  Since $p\in+\mathrm{v}_{n_{0}}$ was arbitrary we
must then have $u_{yy}|_{+\mathrm{v}_{n_{0}}}=0$.  Moreover, the
definitions of $a,b,c$, and $f$ in the regions $-V_{n_{0}}$, $\pm
H_{n_{0}}$ is such that the same arguments used in sections $\S 1$
and $\S 2$ may be applied to show that
\begin{equation*}
u_{yy}|_{-\mathrm{v}_{n_{0}}}=0,\text{ }\text{ }\text{ }\text{
}\text{ } u_{xx}|_{\pm\mathrm{h}_{n_{0}}}=0.
\end{equation*}
Lastly, since $X^{n_{0}}$ was chosen arbitrarily inside the domain
of existence of $u$, (0.4) is valid.  The theorem now follows from
the arguments at the end of the introduction.\medskip

\bigskip

\textit{\textbf{Acknowledgment.}}  Part of this work was carried
out at the University of Pennsylvania while working on my
dissertation.  I would like to thank my advisor, Professor Jerry
Kazdan, as well as Professors Dennis DeTurck and Richard Schoen
for their suggestions and assistance.  This research was supported
by an NSF Postdoctoral Fellowship.

\begin{center}\bigskip
\textsc{References}
\end{center}\bigskip

\noindent[1]\hspace{.06in} \textsc{Y. Egorov}, \textit{On an
example of a linear hyperbolic equation without
solutions},\par\hspace{.01in} C. R. Acad. Sci. Paris S\'{e}r. I
Math. $\mathbf{317}$ no. 12 (1993), 1149-1153.\smallskip

\noindent[2]\hspace{.06in} \textsc{Q. Han}, \textit{Local
isometric embedding of surfaces with Gauss curvature
changing}\par\hspace{.01in} \textit{sign stably across a curve,}
to appear in Calc. Var. \& P.D.E.\smallskip

\noindent[3]\hspace{.06in} \textsc{Q. Han \textit{\&} J.-X. Hong
\textit{\&} C.-S. Lin}, \textit{Local isometric embedding of
surfaces}\par\hspace{.01in} \textit{with nonpositive curvature,}
J. Differential Geom. $\mathbf{63}$ (2003), 475-520.\smallskip

\noindent[4]\hspace{.06in} \textsc{M. Khuri}, \textit{Local
existence of hypersurfaces in $\mathbb{R}^{3}$ with prescribed
Gaussian}\par\hspace{.06in}\textit{curvature at a nondegenerate
critical point}, Electron. J. Diff. Eqns., to appear.\smallskip

\noindent[5]\hspace{.06in} \textsc{M. Khuri}, \textit{The local
isometric embedding in $\mathbb{R}^{3}$ of two-dimensional
Riemannian}\par\hspace{.005in} \textit{manifolds with Gaussian
curvature changing sign to finite order on a
curve},\par\hspace{.01in} J. Differential Geom., to
appear.\smallskip

\noindent[6]\hspace{.06in} \textsc{H. Lewy}, \textit{An example of
a smooth linear partial differential equation
without}\par\hspace{.005in} \textit{solution}, Ann. of Math.
$\mathbf{66}$ (1957), 155-158.\smallskip

\noindent[7]\hspace{.06in} \textsc{C.-S. Lin}, \textit{The local
isometric embedding in $\mathbb{R}^{3}$ of 2-dimensional
Riemannian}\par\hspace{.01in} \textit{manifolds with nonnegative
curvature}, J. Differential Geom. $\mathbf{21}$ (1985),
213-230.\smallskip

\noindent[8]\hspace{.06in} \textsc{C.-S. Lin}, \textit{The local
isometric embedding in $\mathbb{R}^{3}$ of two-dimensional
Riemannian}\par \textit{ manifolds with Gaussian curvature
changing sign cleanly}, Comm. Pure Appl.\par\text{ }Math.
$\mathbf{39}$ (1986), 867-887.\smallskip

\noindent[9]\hspace{.08in} \textsc{N. Nadirashvili}, \textit{The
local embedding problem for surfaces}, preprint.\smallskip

\noindent[10]  \textsc{N. Nadirashvili \textit{\&} Y. Yuan},
\textit{Counterexamples for local isometric
embedding},\par\hspace{.005in} preprint.\smallskip

\noindent[11]  \textsc{O. A. Oleinik}, \textit{On the Cauchy
problem for weakly hyperbolic equations}, Comm.\par\hspace{.01in}
Pure Appl. Math. $\mathbf{23}$ (1970), 569-586.\bigskip

\bigskip
\bigskip
\begin{small}
\noindent\textsc{Marcus A. Khuri:}\newline Department of
Mathematics\newline Stanford University\newline 450 Serra Mall,
Bldg. 380\newline Stanford, CA. 94305-2125\newline
\textsc{E-mail:} \texttt{khuri@math.stanford.edu}
\end{small}

\end{document}